\documentclass[12pt,reqno]{amsart}
\usepackage[dvips]{epsfig}
\usepackage[usenames]{color}
\usepackage{graphicx}
\usepackage{amssymb}
\usepackage{amsthm}
\usepackage{amsmath,epsf}

\numberwithin{equation}{section}
\newtheorem{theo}{Theorem}[section]
\newtheorem{prop}[theo]{Proposition}

\newtheorem{lemm}[theo]{Lemma}

\newtheorem{rema}[theo]{Remark}

\newcommand{\maj}{\mathrm{maj}}  
  
\def\charac{{\raise 2pt\hbox{$\chi$}}}

\def\H{{\mathcal H}}
\def\I{{\mathcal I}}
\def\J{{\mathcal J}}

\def\L{{\mathcal L}}

\def\S{{\mathcal S}}
\def\F{{\mathcal F}}

\def\C{{\mathbb C}}
\def\N{{\mathbb N}}

\def\Young#1{\vbox{\smallskip\offinterlineskip
    \halign{&\vbox{##}\kern-\Thickness\cr #1}}}

\newdimen\Squaresize \Squaresize=10pt
\newdimen\Thickness \Thickness=.3pt
\newdimen\Correction \Correction=7pt

\def\Vide#1{\hbox{
       \vbox to \Squaresize{\vss
          \hbox to \Squaresize{\hss#1 \hss}\vss}
    \hskip-\Correction}
   \kern-\Thickness}

\def\Carre#1{\hbox{\vrule width \Thickness
   \vbox to \Squaresize{\hrule height \Thickness\vss
      \hbox to \Squaresize{\hss#1\hss}
   \vss\hrule height\Thickness}
   \unskip\vrule width \Thickness}
   \kern-\Thickness}

\def\rouge{\textcolor{BrickRed}}
\def\bleu{\textcolor{NavyBlue}}
\def\bleuf{\textcolor{Violet}}
\def\bleup{\textcolor{Periwinkle}}
\def\rose{\textcolor{Orange}}
\def\black{\textcolor{Black}}
\def\violet{\textcolor{Purple}}

\def\jaune{\textcolor{Yellow}}
\def\carre{\bleu{\linethickness{4mm}\line(1,0){.9}}}
\def\rb#1{{\bf\jaune{\footnotesize#1}}}
\def\gris{\textcolor{Gray}}

\title[Coinvariants of $S_n\times S_n$]{Decomposition of the Diagonal Action of $S_n$ on the \\ Coinvariant Space of $S_n\times S_n$}
\author{F. Bergeron}
\address[F. Bergeron]{D\'epartement de Math\'ematiques\\ Universit\'e
  du Qu\'ebec \`a Montr\'eal\\ Montr\'eal, Qu\'ebec, H3C 3P8, CANADA}
\email{bergeron.francois@uqam.ca}
\author{ F. Lamontagne}
\address[F. Lamontagne]{D\'epartement de Math\'ematiques\\ Universit\'e
  du Qu\'ebec \`a Montr\'eal\\ Montr\'eal, Qu\'ebec, H3C 3P8, CANADA}
\email{lamontag@math.uqam.ca}
\date{\today}
\thanks{F. Bergeron is supported in part by NSERC-Canada and
  FQRNT-Qu\'ebec.}

\begin{document}
\begin{abstract}The purpose of this paper is to give an explicit description of the irreducible decomposition of the multigraded $S_n$-module of coinvariants of $S_n\times S_n$.  Many of the results presented can be extended to $S_n^k$, and analogous questions for other finite reflection group.
\end{abstract}
\maketitle
 \parskip=0pt
{\footnotesize \tableofcontents}
\parskip=8pt  

\section{Introduction.} \label{intro}
We study, in this paper, the diagonal action of $S_n$ on the space of coinvariants of $S_n^k$, with main emphasize on the case $k=2$.
We first recall some results of Artin, Shephard-Todd, and Steinberg \cite{artin, shephard, steinberg} that hold for any finite reflection group, in order to specialize them to the symmetric group $S_n$, as well as to $S_n\times S_n$. 

Let $G$ be a finite subgroup of $GL(V)$, where $V$ is a finite vector space with given basis $\mathbf{x}=x_1,\ldots,x_n$. The $x_i$'s are here considered as variables. The group $G$ acts naturally on polynomials $P(\mathbf{x})$ in $\C[\mathbf{x}]$ as
      $$g\cdot P(\mathbf{x}) :=P(g \mathbf{x}),$$
The space of $G$-invariant polynomials, for this action, is denoted $\C[\mathbf{x}]^G$. We also denote $\I_G$ the ideal (of $\C[\mathbf{x}]$) generated by constant term free $G$-invariant polynomials.
By definition, the {\sl coinvariant space} of  $G$ is
\begin{equation}\label{coinvariants}
    \C[\mathbf{x}]_G:=\C[\mathbf{x}]/\I_G.
\end{equation}
It is naturally graded with respect to degree and comes equipped with an action of $G$, since $\I_G$ is both invariant  and homogeneous. A theorem of Steinberg states the group $G$ is generated by reflections if and only if its coinvariant space is isomorphic to the (left) regular representation of $G$.

In particular, for $G=S_n$ the group of permutations, one can  show that $\I_n\ (=\I_{S_n})$ admits the set
      $$\{ h_k(x_k,\ldots, x_n)\ | \ 1\leq k\leq n\ \} $$
 as a Gr\"obner basis. Here we are considering the lexicographic order on monomials, with the usual order $x_1<\cdots <x_n$
 on variables. It follows from general theory  that $ \C[\mathbf{x}]_{S_n}$ can be identified\footnote{In fact it is a set of representatives.} to the linear span of the set
    $$\{ \mathbf{x}^\mathbf{a}\ |\ \mathbf{a}=(a_1,\ldots,a_n)\in \N^n,\ 0\leq a_i<i\ \},$$
with the classical vectorial notation for monomials: 
  $$\mathbf{x}^\mathbf{a}=x_1^{a_1}\ldots x_n^{a_n}.$$
This is often called the {\sl Artin basis} of   $ \C[\mathbf{x}]_{S_n}$. It is evidently an homogeneous ``basis''. This basis make it apparent that the coinvariant space of $S_n$ has dimension $n!$. 

Another important property, of coinvariant spaces of finite groups generated by reflections, is that there is an equivalence of graded $G$-modules
\begin{equation}
\label{tensoriel}
   \C[\mathbf{x}]\simeq  \C[\mathbf{x}]^G\otimes  \C[\mathbf{x}]_G.
\end{equation}
In other words, $\C[\mathbf{x}]$ is a free $ \C[\mathbf{x}]^G$-module. Thus every polynomial $P(\mathbf{x})$ can be expressed in an unique manner (say with respect to the descent basis)  as
\begin{equation}
\label{decomposition}
    P(\mathbf{x})=\sum_{\sigma\in S_n}  f_\sigma(\mathbf{x})\,\mathbf{x}_\sigma,
\end{equation}
with the coefficients $f_\sigma(\mathbf{x})$ a symmetric polynomial.

We recall that the Hilbert series of $\C[\mathbf{x}]$ is easily shown to be
       $$\dim_q \C[\mathbf{x}]=\frac{1}{(1-q)^n},$$
 and that is well known that
      $$\dim_q \C[\mathbf{x}]^{\S_n} =\frac{1}{(q;q)_n},$$
where
     $$(q;q)_n:=(1-q)(1-q^2)\cdots (1-q^n).$$      
It follows from identity (\ref{tensoriel}), in the particular case of $S_n$, that 
  \begin{equation}\label{hilbert_Hn}
     \dim_{q} \C[\mathbf{x}]_{S_n} =   \frac{(q;q)_n}{(1-q)^n}.
\end{equation}
This is readily generalized to other groups generated by reflections.

For all $G$, the coinvariant space $\C[\mathbf{x}]_G$ is always isomorphic (as a graded $G$-module) to the space $\H_G=\I_G^\perp$ of {\sl harmonic polynomials} for $G$. Here $\I_G^\perp$ denotes the orthogonal complement of $\I_G$ for the scalar product on $\C[\mathbf{x}]$ defined as:
\begin{equation}
    \langle P,Q\rangle=P(\partial \mathbf{x})Q(\mathbf{x})\big|_{\mathbf{x}=0}
\end{equation}
In this last expression,  $P(\partial \mathbf{x})$ is to be understood as the differential operator obtained by replacing each variable $x_i$ in $P(\mathbf{x})$ by the partial derivative 
$\partial x_i$, with respect to the variable $x_i$. Moreover, $\mathbf{x}=0$ stands for
    $$x_1=\cdots=x_n=0.$$
 From the fact that $\I_G$ is an ideal, it follows easily that $P(\mathbf{x})$ is in $\H_G$ if and only if 
 it satisfies the system of differential equations
      $$f_a(\partial \mathbf{x}) P(\mathbf{x})=0,\qquad a\in A,$$
where $\{f_a\}_{a\in A}$ is any set of generators for $\C[\mathbf{x}]^G$. In a way, this is equivalent to the fact that the space $\H_G$ is closed under partial derivatives. This leads us to yet another characterization of groups generated by reflections. Namely, there exists an explicit polynomial $\Delta_G(\mathbf{x})$ (see \cite{orbit}), such that the set of all partial derivatives of $\Delta_G$ (of all orders) contains a basis of $\H_G$ if and only if $G$ is a group generated by reflections.

In particular, we can explicitly describe $\H_n\ (=\H_{S_n})$ as the
linear span of all partial derivatives of the Vandermonde determinant $\Delta_n(\mathbf{x})$,
where as usual
   $$\Delta_n(\mathbf{x})=\prod_{i<j} (x_i-x_j).$$
In formula,
\begin{equation}
\label{Ldel}
   \H_n=\L_{\partial} [\Delta_n(\mathbf{x})]
\end{equation}
where the $\L_{\partial} $ in the right hand side stands for ``linear span of all derivatives of''. 

The context for this paper is a ``bigraded'' version of the constructions outlined above  associated to $S_n\times S_n$. More precisely, for $\mathbf{y}=y_1,\dots,y_n$ a second set of $n$ variables, we consider the ring $R:=\C[\mathbf{x},\mathbf{y}]$
of polynomials in the variables $\mathbf{x}$ and $\mathbf{y}$. Here the group $S_n\times  S_n$ acts as a reflection group on $R$ by permuting these two sets of variables independently. Namely,   $$(\sigma,\tau)\, x_i=x_{\sigma(i)},\qquad {\rm and}\qquad (\sigma,\tau)\,y_i=y_{\tau(i)}.$$
Clearly this action respects the ``bidegree'', where the {\sl bidegree}  of a monomial $\mathbf{x}^\mathbf{a} \mathbf{y}^\mathbf{b}$ is $(|\mathbf{a}|,|\mathbf{b}|)$ with
   $$|\mathbf{a}|=a_1+\ldots+a_n,\qquad |\mathbf{b}|=b_1+\ldots+b_n.$$
Our purpose  is to give an explicit description of the irreducible decomposition of  the coinvariant module $R_{S_n\times S_n}$ (or equivalently the module $\H_{S_n\times S_n}$ of $(S_n\times S_n)$-harmonics) considered as an $S_n$-module under the diagonal action. We want this decomposition to take into account the natural bigrading with respect to bidegree (degree in $\mathbf{x}$ and degree in $\mathbf{y}$). We further plan to make explicit (part of\footnote{This will be completed in an upcoming paper. See section \ref{final}.}) the isomorphism of $S_n$-module corresponding to (\ref{tensoriel}). 
Among other things, such a  decomposition gives rise to many beautiful bijections. One of these is between $n$ element subsets of $\N\times\N$ and triples $(D_\sigma,\lambda,\mu)$, with $\lambda$ and $\mu$ partitions  with at most $n$ parts, and $D_\sigma$ varying in a special class of $n$ element subsets of $\N\times \N$, that we call {\em compact diagrams}. These diagrams are naturally indexed by permutations.

\section{Diagonal context.} \label{diagonal}

To emphasize the constructions outlined in section \ref{intro}  in the case of the $S_n$-module of $(S_n\times S_n)$-coinvariants, consider the bigraded component $R_{j,k}$ of bidegree $(j,k)$ of the ring $R=\C[\mathbf{x},\mathbf{y}]$. Clearly this bigraded component affords as a basis the set of monomials
   $$\{\ \mathbf{x}^\mathbf{a} \mathbf{y}^\mathbf{b}\ |\quad |\mathbf{a}|=j\quad {\rm and}\quad |\mathbf{b}|=k \ \},$$ 
with $\mathbf{a},\mathbf{b}\in \N^n$. In turn, the monomial $\mathbf{x}^\mathbf{a} \mathbf{y}^\mathbf{b}$ is in bijection with the {\em bipartite composition}  
   $$(\mathbf{a},\mathbf{b}):=((a_1,b_1),(a_2,b_2),\ldots,(a_n,b_n)),$$
with some of the $(a_i,b_i)$'s possibly equal to $(0,0)$.
We will sometimes use a $2\times n$ matrix notation (also called {\em two line notation}|) for bipartite composition:
   $$(\mathbf{a},\mathbf{b})=\begin{pmatrix}
                   a_1&a_2&\ldots &a_n\\
                   b_1&b_2&\ldots &b_n\end{pmatrix}  $$
When $|\mathbf{a}|=j$ and $|\mathbf{b}|=k$, we say that $(\mathbf{a},\mathbf{b})$ is a  bipartite composition of $(j,k)$. 
Clearly,
   $$(j,k)=(a_1,b_1)+(a_2,b_2)+\ldots+(a_n,b_n).$$
The dimension of $R_{j,k}$ is thus the number of bipartite compositions of $(j,k)$:
  $$\dim R_{j,k}= {n+j-1 \choose j} {n+k-1\choose k}.$$
The ideal generated by constant term free $(S_n\times S_n)$-invariant polynomials, is bihomogeneous with respect to bidegree. The following is clearly a generator set for this ideal:
    $$\{h_1(\mathbf{x}),h_2(\mathbf{x}),\ldots, h_n(\mathbf{x}),h_1(\mathbf{y}),h_2(\mathbf{x}),\ldots,h_n(\mathbf{y})\}.$$
Using a Gr\"obner basis computation, with the order 
   $$x_1<y_1<x_2<y_2<\ldots <x_n<y_n$$
 on the variables and the lexicographic order on monomials, we can identify the space  $R_{\S_n\times S_n}$ with the linear span of the monomials
      $$\{ \mathbf{x}^\mathbf{a} \mathbf{y}^\mathbf{b}\ |\  a_i<i,\ {\rm and}\ b_j<  j\ \}.$$
This is thus a bihomogeneous basis of $R_{\S_n\times S_n}$. We also call this the Artin basis.

We want to make explicit the $S_n$-modules bigraded isomorphisms
\begin{eqnarray}
\label{tens_carre}
    R&\simeq&R^{\S_n\times S_n}\otimes R_{S_n\times S_n}\\
       &\simeq& R^{\S_n\times S_n}\otimes \H_{S_n\times S_n}
\end{eqnarray}
Observe that $R^{S_n\times S_n}$ is in fact isomorphic to $\Lambda(\mathbf{x})\otimes\Lambda(\mathbf{y})$, where we simply write $\Lambda(\mathbf{x})$ for $\C[\mathbf{x}]^{S_n}$. In other words, for each choice of bihomogeneous basis $\mathcal{B}$ of $R_{S_n\times S_n}$ (or $\H_{S_n\times S_n}$), there is a unique decomposition of polynomials $P(\mathbf{x},\mathbf{y})$ of the form
\begin{equation}\label{decomp_tens}
    P(\mathbf{x},\mathbf{y})=\sum_{b\in \mathcal{B}} f_b\ b(\mathbf{x},\mathbf{y}),
 \end{equation}
with coefficients $f_b=f_b(\mathbf{x},\mathbf{y})$ in $\Lambda(\mathbf{x})\otimes\Lambda(\mathbf{y})$. If the polynomial $P(\mathbf{x},\mathbf{y})$ is bihomogeneous of bidegree $(s,t)$ and $b(\mathbf{x},\mathbf{y})$ is bihomogeneous of bidegree $(u,v)$,  then the $f_b$'s can be expressed in the form
\begin{equation}\label{sym_tens}
    f_b(\mathbf{x},\mathbf{y})=\sum_{\lambda\vdash k\atop  \mu\vdash j} a_{\lambda,\mu} m_\lambda(\mathbf{x}) m_\mu(\mathbf{y}),
\end{equation}
with $k-s-u$ and $j=t-v$.
Here, $\lambda\vdash k$ means that $\lambda=(\lambda_1,\lambda_2,\ldots,\lambda_n)$ is a {\em partition} of $k$.  Recall that this is to say that 
  $$k=|\lambda|:=\lambda_1+\lambda_2+\ldots +\lambda_n,$$
 with the $\lambda_i$'s non negative integers such that 
     $$\lambda_1\geq \lambda_2\geq \ldots\geq \lambda_n\geq 0.$$
 The {\em length} $\ell(\lambda)$  is the number of non zero {\em parts} $\lambda_i$ of $\lambda$.
We underline that the number of variables, $n$, is implicitly involved in this description as an upper bound for  both $\ell(\lambda)$ and $\ell(\mu)$.  
We further recall that  the various basis of symmetric functions  are well known (See \cite{macdonald}) to be naturally indexed by partitions. In particular, one such basis corresponds to the {\em monomial} symmetric polynomial $m_\lambda(\mathbf{x})$ defined as:
   $$m_\lambda(\mathbf{x})=\sum \mathbf{x}^\mathbf{a},$$
 where the sum is over all distinct permutations $\mathbf{a}=(a_1,a_2,\ldots,a_n)$ of $\lambda$.  With this in mind, (\ref{sym_tens}) is simply making explicit the usual description of $\Lambda(\mathbf{x})\otimes\Lambda(\mathbf{y})$ for a particular choice of basis of $\Lambda(\mathbf{x})$ and $\Lambda(\mathbf{y})$.

 For $n=2$, we can illustrate all this as follows. Every bihomogeneous polynomial $P(\mathbf{x},\mathbf{y})$ of bidegree $(j,k)$ can be expressed in a unique manner as
     $$P(\mathbf{x},\mathbf{y})=f_{00}+ f_{10} x_2+ f_{01} y_2 + f_{11} x_2 y_2,$$
 with the $f_{uv}(\mathbf{x},\mathbf{y})$ in $\Lambda_{j-u}(x_1,y_2)\otimes \Lambda_{k-v}(y_1,y_2)$.  Here, $\Lambda_j(\mathbf{x})$ denotes the homogeneous component of degree $j$ of $\Lambda(\mathbf{x})$.
For example, we have
\begin{equation*}
x_1y_2+y_1x_2=m_{{1}} ( \mathbf{y} )\,x_2  + m_{{1}} ( \mathbf{x} )\,y_2
   -2\,x_2y_2
\end{equation*}
\section{Algebraic motivation.} \label{motivation}

Our main reason to study $\H_{S_n\times S_n}$ is the following.
Let $\J_n$ be the ideal generated by all (constant term free) diagonally symmetric (invariant under the diagonal action) polynomials. It is clearly bihomogeneous.
The bigraded $S_n$-module $\mathcal{DH}_n$, obtained as the orthogonal complement of $\J_n$, is called the space of {\sl diagonal harmonics} of $S_n$. The space  $\mathcal{DH}_n$  does not appear as a special case of  harmonics for reflection groups (because of the diagonal aspect), and its structure seems much more complicate.
It has been extensively studied recently (see \cite{lattice, Berg2, gordon, gang, haiman, vanishing}), and has many nice properties. In particular its dimension is $(n+1)^{n-1}$. It is clear that we have the bigraded $S_n$-module inclusion
\begin{equation}
   \mathcal{DH}_n\subset \H_{S_n\times S_n},
 \end{equation}
since  $\J_n$ clearly contains $\Lambda(\mathbf{x})\otimes \Lambda(\mathbf{y})$. Thus a more detail understanding of $\H_{S_n\times S_n}$ will shed light on the structure of $\mathcal{DH}_n$. The point here is that
$\H_{S_n\times S_n}$ (or equivalently $R_{S_n\times S_n}$) is much easier to study then $\mathcal{DH}_n$.

It is easy to see that 
  $$\H_{S_n\times S_n}\simeq \H_n\otimes \H_n$$
  as $S_n$-modules, with the diagonal action of $S_n$ on the right hand side.
From the explicit description (\ref{Ldel}) of $\H_n$, it follows that
   $$\H_{S_n\times S_n}=\L_{\partial} [\Delta_n(\mathbf{x})\,\Delta_n(\mathbf{y})],$$
and the above discussion shows that the (bigraded) Hilbert series of $\H_{S_n\times S_n}$ is
  \begin{equation}\label{hilbert_Hndeux}
     \dim_{q,t} \H_{\S_n\times \S_n} =   \frac{(q;q)_n}{(1-q)^n}\frac{(t;t)_n}{(1-t)^n},
 \end{equation}
 which is a bigraded analog of $n!^2$. Much of this can also be formulate in the context of the coinvariant space $R_{S_n\times S_n}$. It is sometimes easier to work in this later context.

\section{Bigrading and bigraded Frobenius characteristic.}\label{frobenius}
In this section, the term symmetric ``function'' is used to underline that we will be using infinitely many variables. The actual variables, $\mathbf{z}=z_1,z_2,z_3,\ldots$, only play a formal role, and will often omit to mention them. Hence we will denote $s_\lambda$, rather then $s_\lambda(\mathbf{z})$, the usual Schur function. These ``formal'' symmetric functions allow a translation of computations on characters of $S_n$ into the more convenient and effective context of symmetric functions through the {\em Frobenius characteristic map}.

 The $S_n$-modules of the previous section (all submodules of $R$) are {\sl bihomogeneous} with respect to bidegree.  We can thus decompose them as direct sums of their bihomogeneous components, which are obtained using the usual linear projections defined on monomials as
     $$\pi_{j,k}(\mathbf{x}^\mathbf{a} \mathbf{y}^\mathbf{b}):=\begin{cases}
     \mathbf{x}^\mathbf{a} \mathbf{y}^\mathbf{b} & |\mathbf{a}|=j,\ {\rm and} \ |\mathbf{b}|=k \\ \\
    0  & \text{otherwise}.
\end{cases}$$
Recall that the {bigraded Frobenius characteristic} of any invariant bihomogeneous submodule $\mathcal{V}$ of $Q[\mathbf{x},\mathbf{y}]$ is defined to be the symmetric function
\begin{equation}
\mathcal{F}_\mathcal{V}(\mathbf{z};q,t):=\sum_{j,k} q^j t^k
                               \frac{1}{n!} \sum_{\sigma\in S_n}
                               \charac_{\mathcal{V}_{j,k} }(\sigma)\,p_{\lambda(\sigma)},
\end{equation}
where $ \charac_{\mathcal{V}_{j,k} }$ is the character of the bihomogeneous component  $\mathcal{V}_{j,k}=\pi_{j,k}(\mathcal{V})$ of $\mathcal{V}$, and where $\lambda(\sigma)$ denotes the (integer)
partition describing the cycle structure of the permutation $\sigma$. As usual, we have denoted here by       
    $$p_\lambda=p_\lambda(\mathbf{z}):=p_{\lambda_1}(\mathbf{z}) p_{\lambda_2}(\mathbf{z}) \cdots p_{\lambda_k}(\mathbf{z})$$ 
the usual {\sl power sum} symmetric functions, with
  $$p_i(\mathbf{z})=z_1^i+z_2^i+z_3^i+\ldots$$

Irreducible representations of $S_n$ are also indexed by partitions of $n$, and there is a natural indexing of them such that the corresponding Frobenius characteristics are the Schur functions $s_\lambda$.
Thus, when expressed in term of Schur functions, the bigraded Frobenius characteristic of an $S_n$-module $\mathcal{V}$ describes the decomposition into irreducibles of each bihomogeneous component of $\mathcal{V}$. Namely
\begin{equation}
   \mathcal{F}_\mathcal{V}(\mathbf{z};q,t)=\sum_{\mu\vdash n}
                       \sum_{j,k} m^\mu_{j,k}q^j t^k  s_\mu,
\end{equation}
where $m^\mu_{j,k}$ is the multiplicity of the irreducible representation of $S_n$
naturally indexed by $\mu$. In other words, if we denote $\mathcal{V}^\mu$ the {\em isotypic component} of $\mathcal{V}$ of type $\mu$, for $\mu$ a partition of $n$, we get
 \begin{equation}
   \mathcal{F}_\mathcal{V}(\mathbf{z};q,t)=\sum_{\mu\vdash n}
                      f_{\mathcal{V}^\mu}(q,t) s_\mu,
\end{equation}  
where $ f_{\mathcal{V}^\mu}(q,t) $ is the Hilbert series of  $\mathcal{V}^\mu$. Recall that $\mathcal{V}^\mu$ is the submodule of $\mathcal{V}$ made of all its irreducible submodules that have same character indexed by $\mu$.
We will give below explicit expressions for the bigraded Frobenius characteristic of both $\H_{\S_n\times \S_n}$ and $R_{S_n\times S_n}$, using
 {\sl plethystic notations}.

\section{Plethystic notation.}\label{plethysm}
For any symmetric function $f$, expressed in the basis of power sum symmetric functions as:
    $$g:=\sum_\mu g_\mu\,p_\mu,$$
we set
     $$g\left[ \frac{\mathbf{z}}{1-q}\right]:=\sum_\mu g_\mu\,p_\mu\left[ \frac{\mathbf{z}}{1-q}\right],$$
where, for a partition $\mu=\mu_1\mu_2 \cdots \mu_r$,
     $$p_\mu\left[ \frac{\mathbf{z}}{1-q}\right]:=\prod_{i=1}^r \frac{p_{\mu_i} }{1-q^{\mu_i}}.$$
The idea here is to consider a power sum $p_k$ as an operator that raises all variables (between brackets) to the power $k$. Hence, if one thinks that $\mathbf{z}=z_1+z_2+\ldots$, then
    $$p_k(\mathbf{z})=p_k(z_1,z_2,\ldots).$$
Thus, the $\mathbf{z}$ that appears in the above expressions stands for the variables of the symmetric functions that we are dealing with. In fact, in this context, we consider the number of variables to be infinite. 

Now, let $h_n=s_{(n)}$ denote the complete homogeneous symmetric functions, then
\begin{equation}\label{frob_pol}
      \mathcal{F}_{\C[\mathbf{x}]} (\mathbf{z};q) = h_n  \left[ \frac{\mathbf{z}}{1-q}\right] .
\end{equation}
is the graded Frobenius characteristic of the $S_n$-module of polynomials.
Hence from \ref{tensoriel}, one deduces  that the graded Frobenius characteristic of $\C[\mathbf{x}]_{\S_n}$ is
\begin{eqnarray}\label{frob_harm}
    \mathcal{F}_{\C[\mathbf{x}]_{\S_n}} (\mathbf{z};q)&=&(q;q)_n h_n  \left[ \frac{\mathbf{z}}{1-q}\right]\\  
          &=& \sum_{\lambda\vdash n} f_\lambda(q) s_\lambda.
 \end{eqnarray}
Thus the coefficients $f_\lambda(q)$ are polynomials with positive integer coefficients\footnote{In fact, the left hand side of (\ref{frob_harm}) is a special case of a Hall-Littlewood symmetric function (see Macdonald \cite{macdonald} for more details).}.

Moreover, observe that
  \begin{equation}
\label{Hall_littlewood}
  \lim_{q\rightarrow 1} \  (q;q)_n h_n  \left[ \frac{\mathbf{z}}{1-q}\right] =h_1^n,
\end{equation}
which is well known to be the Frobenius characteristic of the regular representation of $S_n$. It follows that each  $f_\lambda(q)$ specializes, at $q=1$,  to the  multiplicity $f_\lambda$ of the irreducible indexed by $\lambda$ in the regular representation. It is well known that these $f_\lambda$'s are given by the {\sl hook length formula}.

\section{$S_n$-Frobenius characteristic of the module of ${S_n\times S_n}$ coinvariants.}\label{SnxSn}
One deduces, from results mentioned in section \ref{intro}, that the (simply) graded Frobenius characteristic of $\H_{\S_n}$ is precisely
\begin{equation}\label{harmonique}
    \mathcal{F}_{\H_{\S_n}} (\mathbf{z};q)=(q;q)_n h_n  \left[ \frac{\mathbf{z}}{1-q}\right]
 \end{equation}
In view of (\ref{tens_carre}), we can use this to calculate the bigraded Frobenius characteristic of $\H_{\S_n\times \S_n}$ using the following well known fact. For two $S_n$-modules $\mathcal{V}$ and $\mathcal{W}$,
\begin{equation}
\label{tens_interne}
   \mathcal{F}_{\mathcal{V}\otimes \mathcal{W}}= \mathcal{F}_{\mathcal{V}} *  \mathcal{F}_{ \mathcal{W}},
\end{equation}
where ``$*$'' stands for the ``internal product'' of symmetric functions.
Recall that the {\sl internal product}  of two symmetric functions is the bilinear product such that
  $$p_\lambda* p_\mu = \begin{cases} z_\mu \,p_\mu &{\rm if}\   \lambda=\mu,\\ \\
           0 &{\rm otherwise,}
   \end{cases}$$
where
   $$z_\mu=1^{k_1} k_1! 2^{k_2} k_2! \cdots n^{k_n} k_n! ,$$
 if $\mu$ has $k_i$ parts of size $i$.   We thus easily obtain the following expression for the bigraded Frobenius characteristic of $\H_{\S_n\times \S_n}$, since

\begin{theo}\label{theo:frob des harmo de sn croix sn}
We have
\begin{equation}\label{frob_Hn}
     F_n(\mathbf{z};q,t):=\F_{\H_{S_n\times S_n}}(\mathbf{z};q,t)=(q;q)_n(t;t)_n\ h_n\left[ \frac{\mathbf{z}}{(1-t)(1-q)}\right].
\end{equation}
\end{theo}

\noindent For example, we have
\begin{eqnarray*}
F_{{1}}(\mathbf{z};q,t)&=&s_{{1}}\\
F_{{2}}(\mathbf{z};q,t)&=& ( q\,t+1 ) s_{{2}}+ ( q+t ) s_{{11}}\\
F_{{3}}(\mathbf{z};q,t)&=& ( {q}^{3}{t}^{3}+{q}^{2}{t}^{2}+{q}^{2}t+q{t}^{2}+qt+1
 ) s_{{3}} \\ && + ( {q}^{3}{t}^{2}+{q}^{2}{t}^{3}+{q}^{3}t+{q}^{2
}{t}^{2}+q{t}^{3}+{q}^{2}t+q{t}^{2}+{q}^{2}+qt+{t}^{2}+q+t ) s_{
{21}}\\ && + ( {q}^{2}{t}^{2}+{q}^{3}+{q}^{2}t+q{t}^{2}+{t}^{3}+qt
 ) s_{111}
\end{eqnarray*}
Writing the coefficients of these symmetric functions as matrices gives a better idea of the nice symmetries involved in these expressions.
This is to say that the coefficient of $q^it^j$ is the entry in position $(i,j)$, starting at $(0,0)$ and going from bottom to top and left to right.
Using this convention, $F_4(\mathbf{z};q,t)$ equals to
\begin{eqnarray*}\renewcommand{\arraystretch}{.4}\tiny
\left[ \begin {array}{ccccccc} 0&0&0&0&0&0&1\\\noalign{\medskip}0&0&0&1&1&1&0\\
\noalign{\medskip}0&0&1&1&2&1&0\\\noalign{\medskip}0&1&1&2&1
&1&0\\\noalign{\medskip}0&1&2&1&1&0&0\\\noalign{\medskip}0&1&1&1&0&0&0
\\\noalign{\medskip}1&0&0&0&0&0&0\end {array} \right] s_{{4}}+ \left[
\begin {array}{ccccccc} 0&0&0&1&1&1&0\\\noalign{\medskip}0&1&2&2&2&1&1
\\\noalign{\medskip}0&2&3&4&3&2&1\\\noalign{\medskip}1&2&4&4&4&2&1
\\\noalign{\medskip}1&2&3&4&3&2&0\\\noalign{\medskip}1&1&2&2&2&1&0
\\\noalign{\medskip}0&1&1&1&0&0&0\end {array} \right] s_{{31}}+
 \left[ \begin {array}{ccccccc} 0&0&1&0&1&0&0\\\noalign{\medskip}0&1&1
&2&1&1&0\\\noalign{\medskip}1&1&2&2&2&1&1\\\noalign{\medskip}0&2&2&4&2
&2&0\\\noalign{\medskip}1&1&2&2&2&1&1\\\noalign{\medskip}0&1&1&2&1&1&0
\\\noalign{\medskip}0&0&1&0&1&0&0\end {array} \right] s_{{22}}\\ \\ + \renewcommand{\arraystretch}{.5}\tiny
 \left[ \begin {array}{ccccccc} 0&1&1&1&0&0&0\\\noalign{\medskip}1&1&2
&2&2&1&0\\\noalign{\medskip}1&2&3&4&3&2&0\\\noalign{\medskip}1&2&4&4&4
&2&1\\\noalign{\medskip}0&2&3&4&3&2&1\\\noalign{\medskip}0&1&2&2&2&1&1
\\\noalign{\medskip}0&0&0&1&1&1&0\end {array} \right] s_{{211}}+
 \left[ \begin {array}{ccccccc} 1&0&0&0&0&0&0\\\noalign{\medskip}0&1&1
&1&0&0&0\\\noalign{\medskip}0&1&2&1&1&0&0\\\noalign{\medskip}0&1&1&2&1
&1&0\\\noalign{\medskip}0&0&1&1&2&1&0\\\noalign{\medskip}0&0&0&1&1&1&0
\\\noalign{\medskip}0&0&0&0&0&0&1\end {array} \right] s_{{1111}}
\end{eqnarray*}
We can easily reformulate formula (\ref{frob_Hn}) as
\begin{eqnarray}\label{coeff_qt}
\F_{\H_{S_n\times S_n}}(\mathbf{z};q,t) & = &  \sum_{\lambda\vdash n}  f_\lambda(q,t) s_\lambda\nonumber \\
 & = &\sum_{\lambda\vdash n} ( (q;q)_n(t;t)_n\ s_\lambda\left[\frac{1}{(1-q)(1-t)}\right]) s_\lambda,
\end{eqnarray}
Thus we have an explicit expression for the bigraded enumeration, $f_\lambda(q,t) $, of irreducible representations indexed by $\lambda$ in $\H_{\S_n\times \S_n}$. In particular, using results that can be found in Stanley \cite{stanley2}, and the fact that it corresponds to the coefficient of $s_{1^n}$ in (\ref{coeff_qt}), one finds that the bigraded dimensions of the spaces
 $\mathcal{T}_n$, of {\sl diagonally symmetric} polynomials, and 
 $\mathcal{A}_n$, of {\sl diagonally antisymmetric} polynomials in $\H_{\S_n\times \S_n}$, are respectively
\begin{eqnarray}\label{dim_qt_Triv}
   f_{(n)}(q,t)&=&(q;q)_n(t;t)_n\, h_n\left[ \frac{1}{(1-t)(1-q)}\right]\\
                      &=&\sum_{\sigma\in S_n} q^{{\rm maj}(\sigma)} t^{{\rm maj}(\sigma^{-1})}.\nonumber
\end{eqnarray}
 and
\begin{eqnarray}\label{dim_qt_Alt}
   f_{1^n}(q,t)&=&(q;q)_n(t;t)_n\, e_n\left[ \frac{1}{(1-t)(1-q)}\right]\\
                       &=&\sum_{\sigma\in S_n} q^{{\rm maj}(\sigma)} t^{{n\choose 2}-{\rm maj}(\sigma^{-1})}.\nonumber
\end{eqnarray}
Here ${\rm maj}(\sigma)$ denotes the {\sl major index} of $\sigma$ defined as:
\begin{equation}\label{majeur}
   {\rm maj}(\sigma):=\sum_{i\atop \sigma(i)>\sigma(i+1)} i.
\end{equation}
On the other hand, it is well known (see for instance \cite{stanley2}) that
  $$\sum_{\sigma\in S_n} q^{{\rm maj}(\sigma)} =
   \prod_{k=1}^n\frac{q^k-1}{q-1},$$
which is the classical $q$--analog of $n!$.  In general, using (\ref{Hall_littlewood}), one can easily derive
that
\begin{eqnarray}
\label{coefficients_un}
           f_{\lambda}(q,1)& =&f_\lambda \prod_{k=1}^n\frac{q^k-1}{q-1}, \qquad
           f_{\lambda}(1,t)= f_\lambda\prod_{k=1}^n\frac{t^k-1}{t-1},\\ \label{coefficients_zero}
           f_{\lambda}(q,0)& =&f_\lambda , \qquad\  {\rm and} \qquad\ \
           f_{\lambda}(0,t)= f_\lambda.
\end{eqnarray}
Our purpose in the following sections is to give an explicit combinatorial description of the various isotypic components of
$\H_{\S_n\times \S_n}$ (or equivalently of $R_{S_n\times S_n}$), as an $S_n$-module. However, before going on with this task, let us
recall briefly (see \cite{vanishing} for more details) how to compute the bigraded Frobenius characteristic of $ \mathcal{DH}_n$.

\begin{theo}[Haiman 2002] Let $\nabla$ be the linear operator defined in terms of the modified Macdonald symmetric functions $\widetilde{H}(\mathbf{z};q,t)$
by
\begin{equation}
\label{nabla}
    \nabla \widetilde{H}(\mathbf{z};q,t)= t^{n(\mu)} q^{n(\mu')} \widetilde{H}(\mathbf{z};q,t)
\end{equation}
where $n(\mu):=\sum_k (k-1)\,\mu_k$. Then we have
\begin{equation}
\label{frob_Diag}
   \mathcal{F}_{\mathcal{DH}_n} (\mathbf{z};q,t)= \nabla(e_n)
\end{equation}
where $e_n$ is the $n^{\rm th}$ elementary symmetric function.
\end{theo}
The operator $\nabla$, defined by (\ref{nabla}), has been introduced by the first author with A.~Garsia to study properties of Macdonald polynomials in conjunction with the study of diagonal harmonics. Many conjectures about it remain open (See \cite{nabla}).
For example, we have
\begin{eqnarray*}
\nabla(e_{{1}} ) &=&s_{{1}}\\
\nabla(e_{{2}} ) &=&s_{{2}}+ ( q+t ) s_{{11}}\\
\nabla(e_{{3}} ) &=&s_{{3}}+ ( {q}^{2}+qt+{t}^{2}
+q+t ) s_{{21}}+ ( {t}^{3}+q{t}^{2}+{q}^{2}t+qt+{q}^{3} ) s_{{111}}
\end{eqnarray*}
Recently, an explicit conjecture (see \cite{gang}) for a combinatorial description of the coefficients of $\nabla(e_n)$, when expressed in term of the {\sl monomial basis} $m_\mu$, has been suggested. It should be interesting to show how this combinatorial description can be explained in the larger context of  $R_{S_n\times S_n}$. Illustrating with $n=3$, this would make evident why the coefficients of the $m_\lambda$'s in
$$\nabla(e_3)=m_{{3}}+ \renewcommand{\arraystretch}{.5}
\left[ \begin {array}{ccc} 1&0&0\\\noalign{\medskip}1&1&0\\\noalign{\medskip}
1&1&1\end {array} \right] m_{{21}}+ \left[ \begin {array}{cccc} 1&0&0
&0\\\noalign{\medskip}2&1&0&0\\\noalign{\medskip}2&3&1&0
\\\noalign{\medskip}1&2&2&1\end {array} \right] m_{{111}}$$
are ``contained'' in the corresponding coefficients in
  $$F_3(\mathbf{z};q,t)=\renewcommand{\arraystretch}{.5}\left[ \begin {array}{cccc} 0&0&0&1\\\noalign{\medskip}0&1&1&0
\\\noalign{\medskip}0&1&1&0\\\noalign{\medskip}1&0&0&0\end {array}
 \right] m_{{3}}+ \left[ \begin {array}{cccc} 0&1&1&1
\\\noalign{\medskip}1&2&2&1\\\noalign{\medskip}1&2&2&1
\\\noalign{\medskip}1&1&1&0\end {array} \right] m_{{2,1}}+ \left[
\begin {array}{cccc} 1&2&2&1\\\noalign{\medskip}2&4&4&2
\\\noalign{\medskip}2&4&4&2\\\noalign{\medskip}1&2&2&1\end {array}
 \right] m_{{1,1,1}}
$$

\section{Properties of the bigraded Frobenius of coinvariants.}\label{hilbert}
Many symmetries are apparent in the explicit values of bigraded Frobenius Characteristics $F_n(\mathbf{z};q,t)$ given above. In fact, we have the following.

\begin{prop}\label{symetries}
  For all $n$, we have
 \begin{enumerate}
     \item $F_n(\mathbf{z};t,q) = F_n(\mathbf{z};q,t) $
      \item $F_n(\mathbf{z};t,q) = (q\,t)^{n\choose 2} F_n(\mathbf{z};q^{-1},t^{-1})$
      \item $\begin{array}[t]{lcl}
               \omega(F_n(\mathbf{z};q,t)) & = & q^{n\choose 2} F_n(\mathbf{z};q^{-1},t) \\
                                    & = & t^{n\choose 2} F_n(\mathbf{z};q,t^{-1})
                \end{array}$
\end{enumerate}
\end{prop}
Here $\omega$ is the usual involution on symmetric functions that send $s_\lambda$ to $s_{\lambda'}$ ($\lambda'$ denoting, as usual, the conjugate of $\lambda$). This is well known to correspond to ``twisting'' by the sign representation.

{\bf Sketch of Proof.} All of these symmetries correspond to automorphisms (or anti automorphisms) of $\H_{\S_n\times \S_n}$.
\vskip-50pt
 \begin{enumerate}
     \item corresponds to the evident symmetry of $\H_{\S_n\times \S_n}$  that corresponds to exchanging the $\mathbf{x}$ variables with the $\mathbf{y}$ variables.

     \item corresponds to the morphism that sends $P(\mathbf{x},\mathbf{y})$ into $P(\partial \mathbf{x},\partial \mathbf{y})\Delta_n(\mathbf{x})\Delta_n(\mathbf{y})$.

     \item corresponds to the anti automorphism (twisting by the sign representation) that sends $P(\mathbf{x},\mathbf{y})$ into
     $P(\partial \mathbf{x},\mathbf{y}) \Delta_n(\mathbf{x})$, for the first identity.  
     
     \item[3')] For the second identity, we rather consider the anti automorphism
$P(\mathbf{x},\mathbf{y})\mapsto  P(\mathbf{x},\partial \mathbf{y}) \Delta_n(\mathbf{y})$.
\end{enumerate}

We will respectively denote these last two anti automorphism
  \begin{equation}\label{flipX}
       \downarrow_{\mathbf{x}}\!{P(\mathbf{x},\mathbf{y})}:= P(\partial\, \mathbf{x},\mathbf{y}) \Delta_n(\mathbf{x})
  \end{equation}
 and
    \begin{equation}\label{flipY}
      \downarrow_{\mathbf{y}}\!{P(\mathbf{x},\mathbf{y})}:= P(\mathbf{x},\partial\, \mathbf{y}) \Delta_n(\mathbf{y}).
  \end{equation}
These are both called ``flips''; with respect to $\mathbf{x}$ in the first case, and $\mathbf{y}$ in the second. Clearly, if $P(\mathbf{x},\mathbf{y})$ is of bidegree $(j,k)$, then   $\downarrow_{\mathbf{x}}\!{P(\mathbf{x},\mathbf{y})}$ is of bidegree $({n\choose 2}-j,k)$, and  $\downarrow_{\mathbf{y}}\!{P(\mathbf{x},\mathbf{y})}$ is of bidegree $(j,{n\choose 2}-k)$.

\section{Trivial component.}\label{isotypic_trivial}
Let us now recapitulate. According to  (\ref{tens_carre}), we have a decomposition of $R$, as a tensor product of the ring of $(S_n\times S_n)$-invariants with the coinvariant module $R_{S_n\times S_n}$, giving a $S_n$-module isomorphism
\begin{equation}\label{decomposition_module}
    R\simeq \Lambda(\mathbf{x})\otimes \Lambda(\mathbf{y})\otimes R_{S_n\times S_n},
 \end{equation}
 for the diagonal action of $S_n$.
Formula (\ref{frob_Hn}) translates this decomposition in term of Frobenius characteristic. Since, both $\Lambda(\mathbf{x})$ and $\Lambda(\mathbf{y})$  are $S_n$ invariant by definition, it follows that each isotypic component $R^\lambda$ of $R$ decomposes as
  \begin{equation}\label{isotypic_decomposition}
    R^\lambda \simeq \Lambda(\mathbf{x})\otimes \Lambda(\mathbf{y})\otimes R_{S_n\times S_n}^\lambda.
  \end{equation}
Although we will now concentrate on the {\em trivial} ($\lambda=(n)$) and {\em alternating} ($\lambda=1^n$) components, most of constructions and results that follow can be extended to all isotypic components\footnote{See comment in section \ref{final}.}. 

For the purpose of our exposition, we need to recall some definitions. An $n$-cell {\em diagram} (also called {\em bipartite partition}  in the literature)
    $$(\alpha,\beta)=((a_1,b_1),(a_2,b_2),\ldots, (a_n,b_n))$$
is a bipartite composition with  {\em cells}, $(a_i,b_i)$, listed in increasing\footnote{It is convenient for our exposition to follow this convention rather than the usual one which corresponds to the decreasing lexicographic order..} {\em lexicographic} order. Recall that this is the order such that   
  $$(a,b)\preceq (a',b')\quad {\rm iff}\quad
   \begin{cases}b<b' & {\rm or,}\\
                b=b'\quad {\rm and}\quad a\leq a'.
   \end{cases}$$
In other words, the two line representation of a diagram
 \begin{equation}\label{two_line}
   (\alpha,\beta)=\begin{pmatrix}
                   a_1&a_2&\ldots &a_n\\
                   b_1&b_2&\ldots &b_n\end{pmatrix} 
 \end{equation}
is such that the columns are in increasing lexicographic order. Observe that a special cases of diagrams corresponds to the usual two line notation for permutations:  
 \begin{equation}\label{graphe}
      \sigma=\begin{pmatrix}
                   1&2&\ldots &n\\
                   \sigma(1)&\sigma(2)&\ldots &\sigma(n)\end{pmatrix} 
 \end{equation}
The trivial isotypic component $R^{S_n}$ of $R$ is clearly spanned by the ($n$-cell diagrams indexed) set of {\em monomial} diagonally symmetric polynomials
     $$M_{(\alpha,\beta)}:=\sum \mathbf{x}^\mathbf{a}\mathbf{y}^\mathbf{b},$$
where the sum is over all distinct bipartite compositions $(\mathbf{a},\mathbf{b})$ obtained by permuting the cells of $(\alpha,\beta)$. For example,
   $$M_{\left({012\atop 001}\right)}= y_2x_3^2y_3+y_1x_3^2y_3+x_2^2y_2x_3+x_1^2y_1x_3+y_1x_2^2y_2
               +y_1x_2x_1^2$$
Observe that the leading monomial of $M_{(\alpha,\beta)}$ is $X^\alpha Y^\beta$, for the
lexicographic monomial order with the variables order:
    $$x_n>y_n>\ldots > x_2>y_2>x_1>y_1.$$
Moreover  $M_{(\alpha,\beta)}$ is clearly bihomogeneous of bidegree $(|\alpha|,|\beta|)$.
For more on diagonally symmetric polynomials, see \cite{gessel} or \cite{rosas}.

In view of (\ref{isotypic_decomposition}),  the trivial isotypic component of the space
     $$ \mathcal{T}_n:=R_{S_n\times S_n}^{S_n}$$
has dimension of $ \mathcal{T}_n$ is $n!$. In fact, its bigraded dimension is given by formula (\ref{dim_qt_Triv}). We will now construct, for each permutation $\sigma$ in $S_n$, a diagram
${(\alpha(\sigma), \beta(\sigma))}$ with the property that
  \begin{equation}\label{maj_prop}
       |\alpha(\sigma)|=\maj(\sigma),\qquad {\rm and}\qquad |\beta(\sigma)|=\maj(\sigma^{-1}).
  \end{equation}
Furthermore, we will verify that any diagonally symmetric polynomial  can be uniquely decomposed in term of the symmetric polynomials associated to these monomials. More explicitly we will show that, for all diagram $(\gamma,\delta)$, we have
  \begin{equation}\label{decomp_triv}
       M_{(\gamma,\delta)}=\sum_{\sigma\in S_n} 
                               f_\sigma(\mathbf{x},\mathbf{y})\, M_{(\alpha(\sigma), \beta(\sigma))},
   \end{equation}
  with $f_\sigma(\mathbf{x},\mathbf{y})$ in $\Lambda(\mathbf{x})\otimes \Lambda(\mathbf{y})$.
 
For example, with $n=3$, our construction will give as a basis for $\mathcal{T}_3$ the set
\begin{eqnarray*}
M_{\left({000\atop 000}\right)}&=&1,\\
M_{\left({011\atop 100}\right)}&=&x_1x_2y_3+x_1y_2x_3+y_1x_2x_3,\\
M_{\left({001\atop 110}\right)}&=&x_1y_2y_3+x_2y_1y_3+x_3y_1y_2,\\
M_{\left({001\atop 010}\right)}&=&x_1y_2+x_1y_3+x_2y_1+x_2 y_3+x_3y_1+x_3y_2, \\
M_{\left({011\atop 101}\right)}&=&
x_1y_1x_2y_3+x_1y_1y_2x_3+x_1x_2y_2y_3+x_1y_2x_3y_3+y_1x_2y_2x_3+y_1x_2x_3y_3,\\
M_{\left({012\atop 210}\right)}&=&
x_1^2x_2y_2y_3^2+x_1^2y_2^2x_3y_3+x_1y_1x_2^2y_3^2+x_1y_1y_2^2x_3^2+y_1^2x_2^2x_3y_3+y_1^2x_2y_2x_3^2.
\end{eqnarray*}
Below are a few illustrations of decompositions of form (\ref{decomp_triv}) using this basis.
\begin{eqnarray*}
M_{\left({001\atop 001}\right)} &=&s_{1} ( \mathbf{x} ) s_{1} ( \mathbf{y} ) \,M_{\left({000\atop 000}\right)}
                    -M_{\left({001\atop 010}\right)} \\
M_{\left({002\atop 001}\right)}  &=&s_{2} ( \mathbf{x} ) s_{1} ( \mathbf{y} )\,M_{\left({000\atop 000}\right)}  
                   -s_{1} ( \mathbf{x} )\,M_{\left({001\atop 010}\right)} 
                   +  M_{\left({011\atop 100}\right)} \\
M_{\left({002\atop 010}\right)}  &=&s_{1} ( \mathbf{x} )\,M_{\left({001\atop 010}\right)}
                    -s_{11} ( \mathbf{x} ) s_{1} ( \mathbf{y} )\,M_{\left({000\atop 000}\right)}
                    -M_{\left({011\atop 100}\right)}  \\
M_{\left({011\atop 001}\right)}   &=&s_{11} ( \mathbf{x} )s_{1} ( \mathbf{y} ) \,M_{\left({000\atop 000}\right)}
                  -M_{\left({011\atop 100}\right)} \\ 
M_{\left({001\atop 002}\right)} &=&s_{1} ( \mathbf{x} ) s_{2} ( \mathbf{y} ) \,M_{\left({000\atop 000}\right)}
                  -s_{1} ( \mathbf{y} )\,M_{\left({001\atop 010}\right)} 
                  +M_{\left({001\atop 110}\right)} \\
M_{\left({001\atop 020}\right)} &=&s_{1} ( \mathbf{y} )\, M_{\left({001\atop 010}\right)}
                  - s_{1} ( \mathbf{x} ) s_{11} ( \mathbf{y} ) \,M_{\left({000\atop 000}\right)}
                  -M_{\left({001\atop 110}\right)}\\
M_{\left({011\atop 001}\right)} &=&s_{1} ( \mathbf{x} ) s_{11} ( \mathbf{y} )  \,M_{\left({000\atop 000}\right)}
                 -M_{\left({001\atop 110}\right)} 
\end{eqnarray*}
For any permutation $\sigma$ in $S_n$, we simply define
\begin{equation}\label{alpha_sigma}
   \alpha(\sigma):=(d_{1}(\sigma),d_{2}(\sigma),\ldots,d_{n}(\sigma)),
\end{equation}
with
 \begin{equation}
     \label{descentes_i}
      d_i(\sigma):= \# \{k\ |\ k<i\quad{\rm and}\quad \sigma_k >\sigma_{k+1}\ \}.
\end{equation}
If we further define
 \begin{equation}
\label{beta_sigma}
   \beta(\sigma):=(d_{\sigma(1)}(\sigma^{-1}),d_{\sigma(2)}(\sigma^{-1}),\ldots,d_{\sigma(n)}(\sigma^{-1})),
\end{equation}
then it is clear that identities (\ref{maj_prop}) hold for the bipartite composition $(\alpha(\sigma),\beta(\sigma))$. The fact that the cells of
$(\alpha(\sigma),\beta(\sigma))$ are actually in increasing lexicographic order is also readily verified. Hence we get a diagram. We will show in section
\ref{base coinvariants} that the set
  \begin{equation}\label{la_base_triviale}
       \mathcal{M}_n:=\{\ M_{(\alpha(\sigma),\beta(\sigma))}\ |\ \sigma\in S_n\ \},
   \end{equation}
is indeed a basis of $\mathcal{T}_n$.

\section{Alternating component.}\label{isotypic_alternating}

We can easily transpose to the submodule $R^{\pm}$ of diagonal alternants of $R$, the discussion of the last section. 
This space affords as linear basis the set of all determinants
  $$\Delta_D (\mathbf{x},\mathbf{y}) := \det \Bigl( (x_i^a y_i^b)_{\substack{(a,b) \in D \\ 1 \leq i \leq n}} \Bigr)$$
with $D$ varying in the set of {\em strict} diagrams. This is to say that we are not allowing repetition of cells. Observe that, choosing $D = \{ (i,0) \ | \ 0 \leq i \leq n-1 \}$, we get the usual Vandermonde determinant as a special case. The module
 of {\em diagonal harmonic alternants}  has been the object of a lot interest (See \cite{gordon, gang, haiman}) in the last 15 years. It is a submodule of
   $$ \mathcal{A}_n = R^{\pm} \cap \mathcal{H}_{S_n \times S_n }. $$
We have already observed that this last module has dimension $n!$ and its Hilbert series is
  $$(q; q)_n (t; t)_n e_n \left[ \frac{1}{(1-q)(1-t)} \right]. $$
The space $R^{\pm}$ is a free $ (\Lambda_n (\mathbf{x}) \otimes \Lambda_n (\mathbf{y}) )$-module, which is made explicit by the isomorphism
  $$R^{\pm} \simeq \Lambda_n (\mathbf{x}) \otimes \Lambda_n (\mathbf{y}) \otimes \mathcal{A}_n. $$
To get a basis of $\mathcal{A}_n$, we need only ``flip'' the basis $\mathcal{M}_n$ of section \ref{isotypic_trivial}, either with respect to $\mathbf{x}$ or $\mathbf{y}$. 

\section{Generalities on diagrams.}\label{compact}
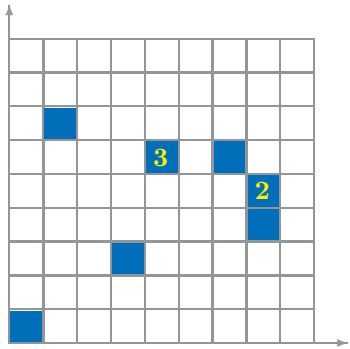
\begin{figure}[ht] 
\gris{
\begin{picture}(10,10)(0,0)
\put(0,0){\vector(1,0){10}}
\put(0,0){\vector(0,1){10}}
\multiput(0,0)(0,1){10}{\line(1,0){9}}
\multiput(0,0)(1,0){10}{\line(0,1){9}}
\put(0.05,0.5){\carre}
\put(1.05,6.5){\carre}
\put(3.05,2.5){\carre}
\put(4.05,5.5){\carre}
\put(4.25,5.25){\rb{3}}
\put(6.05,5.5){\carre}
\put(7.05,3.5){\carre}
\put(7.05,4.5){\carre}
\put(7.25,4.25){\rb{2}}
\end{picture}}
\caption{A diagram.}
\label{exemple_diag}
\end{figure}

We intend to describe a natural classification of $n$-cells diagrams in terms of permutations in $S_n$. To this end, it is sometimes better appropriate to
think of diagrams as $n$-element multisubsets of $\N\times \N$. A diagram $D=(\alpha,\beta)$ can thus be ``drawn'' as a multiset of $1\times 1$ boxes,
in the combinatorial plane $\N\times \N$. Figure \ref{exemple_diag} gives the representation of the diagram whose two line representation is:
\begin{equation}\label{diag_exemple}
    \begin{pmatrix}
       0&1&3&4&4&4&6&7&7&7\\
       0&6&2&5&5&5&5&3&4&4
   \end{pmatrix}
\end{equation}

The numbers appearing in the cells are multiplicities, and when no multiplicity is mentioned it is understood to be 1.

Let $D$ be a $n$-cell diagram
\begin{equation}\label{un_diagramme}
   D=\begin{pmatrix}
                   a_1&a_2&\ldots &a_n\\
                   b_1&b_2&\ldots &b_n\end{pmatrix},
 \end{equation}
 hence its cells are ordered in increasing lexicographic order. 
 We say that $i$ is a {\em descent} of $D$, if $a_{i+1}>a_i$ and $b_{i+1}<b_i$. 
 We denote ${\rm Desc}(D)$ the set of descents of $D$.
For $(a,b)$ in $D$, the set ${\rm Desc}_{(a,b)}(D)$ of descents of $D$ that precede  $(a,b)$ is then defined as
  $${\rm Desc}_{(a,b)}(D):={\rm Desc}(D) \cap \{ k \ |\ (a_k,b_k)\preceq (a,b)\ \}.$$
We denote $d_{(a,b)}(D)$ the cardinality of ${\rm Desc}_{(a,b)}(D)$. Writing simply
$d_i(D)$ for $d_{(a_i,b_i)}(D)$, we see that this definition generalizes our previous notion of definition (\ref{descentes_i}). Observe that, for each $(a,b)$ in $D$, we have
\begin{equation}\label{inegalite}
   a\, \geq\, d_{(a,b)}(D).
\end{equation}

We classify $n$-cell diagrams in term of permutations of $S_n$ in the following way. We associate (as described below)  to each $n$-cell diagram $D$ a certain permutation $\sigma(D)$, and set
\begin{equation}
\label{equivalence}
 D\simeq D',\qquad{\rm iff}\qquad \sigma(D)=\sigma(D').
\end{equation}
We define the {\sl  classifying permutation} $\sigma=\sigma(D)$, of an $n$-cell diagram $D$, to be the unique permutation, $\sigma\in S_n$, such that $\sigma^{-1}$  reorders the $b_i$'s  in increasing order: 
       $$b_{\sigma^{-1}(1)}\leq b_{\sigma^{-1}(2)}\leq \ \ldots\   \leq b_{\sigma^{-1}(n)},$$
 in such a way that $\sigma(i+1)=\sigma(i)+1$,  whenever $b_i=b_{i+1}$ and $a_i\leq a_{i+1}$. Observe that this definition forces the descents of $\sigma$ to be the same as those of $D$. Moreover, the descents of $\sigma^{-1}$ are the same as those of $D^{-1}$.

The classifying permutation of the diagram in (\ref{diag_exemple}) is:
\begin{equation}\label{sa_permutation}
\begin{pmatrix}
 1&2&3&4&5&6&7&8&9&10\\
 1&10&2&6&7&8&9&3&4&5
      \end{pmatrix}
\end{equation}
Evidently a permutation is its own classifying permutation:
    $$\sigma(\tau)=\tau.$$
We further associate to each diagram $D=(\alpha,\beta)$ a special diagram $\Gamma(D)$, called the {\em compactified} of $D$, as follows:
\begin{equation}\label{def_compact}
   \Gamma(D):=\{  (d_{(a,b)}(D), d_{(b,a)}(D^{-1}))\ |\ (a,b)\in D\ \}.
\end{equation}
Here $D^{-1}$ is the diagram obtained by reordering $(\beta,\alpha)$ in increasing lexicographic order. Naturally, $D$ is said to be {\em compact} if and only if $D=\Gamma(D)$. 

\begin{figure}[ht] 
\gris{\begin{picture}(6,6)(0,0)
\put(0,0){\vector(1,0){5}}
\put(0,0){\vector(0,1){5}}
\multiput(0,0)(0,1){5}{\line(1,0){4}}
\multiput(0,0)(1,0){5}{\line(0,1){4}}
\put(0.05,0.5){\carre}
\put(1.05,0.5){\carre}
\put(2.05,0.5){\carre}
\put(2.25,0.25){\rb{3}}
\put(1.05,1.5){\carre}
\put(1.25,1.25){\rb{4}}
\put(0.05,2.5){\carre}
\end{picture}}
\caption{A compact diagram.}
\label{exemple_compact}
\end{figure}
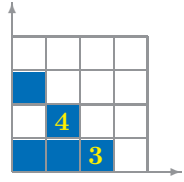

The compactified of the diagram appearing in (\ref{diag_exemple}) is
$$    \begin{pmatrix}
      0&0&1&1&1&1&1&2&2&2\\
      0&2&0&1&1&1&1&0&0&0
   \end{pmatrix}
$$
It corresponds to the pictorial representation of Figure \ref{exemple_compact}. 
The observations above  imply that    
 \begin{equation}\label{characterisation}
       \Gamma(\sigma(D))=\Gamma(D),
    \end{equation}
 since the ``position'' taken by $(b_i,a_i)$ in $D^{-1}$ is $\sigma(i)$.
We will show that the compactified of a diagram is uniquely characterized by its classifying permutation. This implies that there are exactly $n!$ compact $n$-cell diagrams, one in each equivalence class with respect to relation $\simeq$. They can thus be naturally labeled $D_\sigma$, for the corresponding classifying permutation $\sigma$.  In  fact, $D_\sigma$ is none other then the diagram $(\alpha(\sigma),\beta(\sigma))$ considered in the definition of $\mathcal{M}_n$ in section \ref{isotypic_trivial}. Figures \ref{fig:multidiagrams pour 3} and \ref{fig:multidiagrams pour 4} respectively give the compact diagrams for $n$ equal $3$ and $4$, with the corresponding permutation labels. 
\def\carre{\bleu{\linethickness{2.6mm}\line(1,0){.9}}}
\def\rb#1{{\bf\jaune{\tiny#1}}}
\setlength{\unitlength}{3mm}\vskip-15pt
{\tiny \begin{figure}[ht] 
$$\begin{matrix} \qquad
\begin{picture}(1,1)(0,0)
\put(0,0.5){\carre}
\put(0.2,0.25){\rb{3}}
\end{picture} \qquad & \qquad 
\begin{picture}(2,2)(0,0)
\put(0,0.5){\carre}
\put(0,1.5){\carre}
\put(1,0.5){\carre}
\end{picture} \qquad & \qquad
\begin{picture}(2,2)(0,0)
\put(0,1.5){\carre}
\put(1,0.5){\carre}
\put(1,1.5){\carre}
\end{picture} \qquad \\
123 & 132 & 213 \\ 
\begin{picture}(2,2)(0,0)
\put(1,0.5){\carre}
\put(0,1.5){\carre}
\put(0.2,1.25){\rb{2}}
\end{picture} &
\begin{picture}(2,2)(0,0)
\put(0,1.5){\carre}
\put(1,0.5){\carre}
\put(1.2,0.25){\rb{2}}
\end{picture} &
\begin{picture}(3,3.5)(0,0)
\put(0,2.5){\carre}
\put(1,1.5){\carre}
\put(2,0.5){\carre}
\end{picture}\\
231 & 312 & 321
\end{matrix}$$\vskip-10pt
\caption{Compact diagrams for $n=3$.}
\label{fig:multidiagrams pour 3}
\end{figure}
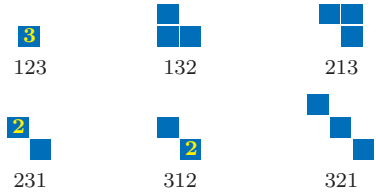}\vskip-10pt
{\tiny
\begin{figure}[ht] 
$$\begin{matrix}
\qquad \begin{picture}(1,1)(0,0)
\put(0,0.5){\carre}
\put(0.1,0.25){\rb{4}}
\end{picture} \qquad & \qquad 
\begin{picture}(2,2)(0,0)
\put(0,0.5){\carre}
\put(0.1,0.25){\rb{2}}
\put(0,1.5){\carre}
\put(1,0.5){\carre}
\end{picture} \qquad & \qquad
\begin{picture}(2,2)(0,0)
\put(0,1.5){\carre}
\put(1,0.5){\carre}
\put(0,0.5){\carre}
\put(1,1.5){\carre}
\end{picture} \qquad & \qquad
\begin{picture}(2,2)(0,0)
\put(1,0.5){\carre}
\put(0,1.5){\carre}
\put(0.1,1.25){\rb{2}}
\put(0,0.5){\carre}
\end{picture} \qquad  \\
  1234 & 1243 & 1324 & 1342 \\ 
\begin{picture}(2,2)(0,0)
\put(0,1.5){\carre}
\put(1,0.5){\carre}
\put(1.1,0.25){\rb{2}}
\put(0,0.5){\carre}
\end{picture} &
\begin{picture}(3,3)(0,0)
\put(0,0.5){\carre}
\put(1,1.5){\carre}
\put(0,2.5){\carre}
\put(2,0.5){\carre}
\end{picture} &
\begin{picture}(2,2)(0,0)
\put(0,1.5){\carre}
\put(1,0.5){\carre}
\put(1,1.5){\carre}
\put(1.1,1.25){\rb{2}}
\end{picture} &
\begin{picture}(3,3.5)(0,0)
\put(0,1.5){\carre}
\put(1,0.5){\carre}
\put(1,2.5){\carre}
\put(2,1.5){\carre}
\end{picture} \\
1423 & 1432 & 2134 & 2143 \\ 
\begin{picture}(2,2)(0,0)
\put(1,0.5){\carre}
\put(0,1.5){\carre}
\put(0.1,1.25){\rb{2}}
\put(1,1.5){\carre}
\end{picture} &
\begin{picture}(2,2)(0,0)
\put(1,0.5){\carre}
\put(0,1.5){\carre}
\put(0.1,1.25){\rb{3}}
\end{picture} &
\begin{picture}(2,3)(0,0)
\put(0,1.5){\carre}
\put(1,0.5){\carre}
\put(0,2.5){\carre}
\put(1,1.5){\carre}
\end{picture} & 
\begin{picture}(3,3.5)(0,0)
\put(0,1.5){\carre}
\put(0,2.5){\carre}
\put(2,0.5){\carre}
\put(1,1.5){\carre}
\end{picture} \\
2314 & 2341 & 2413 & 2431 \\ 
\begin{picture}(2,2)(0,0)
\put(0,1.5){\carre}
\put(1,0.5){\carre}
\put(1.1,0.25){\rb{2}}
\put(1,1.5){\carre}
\end{picture} & 
\begin{picture}(3,2)(0,0)
\put(0,1.5){\carre}
\put(1,0.5){\carre}
\put(2,0.5){\carre}
\put(1,1.5){\carre}
\end{picture} &
\begin{picture}(3,3)(0,0)
\put(2,2.5){\carre}
\put(2,0.5){\carre}
\put(1,1.5){\carre}
\put(0,2.5){\carre}
\end{picture} &
\begin{picture}(3,3.5)(0,0)
\put(0,2.5){\carre}
\put(2,0.5){\carre}
\put(1,1.5){\carre}
\put(1,2.5){\carre}
\end{picture} \\
3124 & 3142 & 3214 & 3241 \\ 
\begin{picture}(2,2)(0,0)
\put(0,1.5){\carre}
\put(0.1,1.25){\rb{2}}
\put(1,0.5){\carre}
\put(1.1,0.25){\rb{2}}
\end{picture} & 
\begin{picture}(3,3)(0,0)
\put(1,1.5){\carre}
\put(0,2.5){\carre}\put(0.1,2.25){\rb{2}}
\put(2,0.5){\carre}
\end{picture} &
\begin{picture}(2,2)(0,0)
\put(0,1.5){\carre}
\put(1,0.5){\carre}
\put(1.1,0.25){\rb{3}}
\end{picture} &
\begin{picture}(3,3.5)(0,0)
\put(1,0.5){\carre}
\put(0,2.5){\carre}
\put(2,0.5){\carre}
\put(1,1.5){\carre}
\end{picture} \\
3412 & 3421 & 4123 & 4132 \\ 
\begin{picture}(3,3)(0,0)
\put(0,2.5){\carre}
\put(2,0.5){\carre}
\put(1,1.5){\carre}
\put(2,1.5){\carre}
\end{picture} &
\begin{picture}(3,3)(0,0)
\put(0,2.5){\carre}
\put(2,0.5){\carre}
\put(1,1.5){\carre}
\put(1.1,1.25){\rb{2}}
\end{picture} &
\begin{picture}(3,3)(0,0)
\put(1,1.5){\carre}
\put(0,2.5){\carre}
\put(2,0.5){\carre}\put(2.1,0.25){\rb{2}}
\end{picture} & 
\begin{picture}(4,4.5)(0,0)
\put(1,2.5){\carre}
\put(2,1.5){\carre}
\put(0,3.5){\carre}
\put(3,0.5){\carre} \end{picture}
\\
4213 & 4231 & 4312 & 4321
\end{matrix}$$\vskip-10pt
\caption{Compact diagrams for $n=4$.}
\label{fig:multidiagrams pour 4}
\end{figure}
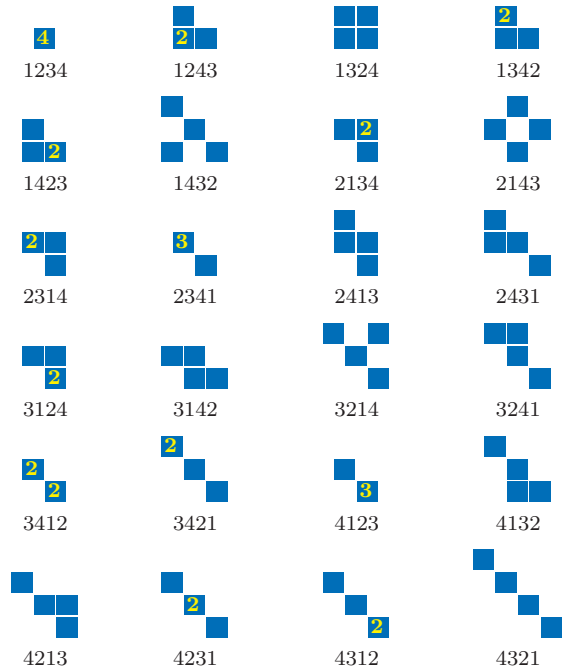}

\section{Characterizations of compact diagrams.}\label{characterization}
We are now going to characterize $D_\sigma$ as being ``minimal'' among those that are equivalent with respect to the relation $\simeq$.  This is the motivation behind the terminology. 

It immediately follows from definition (\ref{def_compact}) that the descent set of $\Gamma(D)$ and of $\Gamma(D)^{-1}$ are respectively the same as those for $D$ and $D^{-a}$. This implies that
 \begin{equation}\label{idemp_gamma}
     \Gamma^2(D)=\Gamma(D).
 \end{equation}
Moreover, in view of (\ref{inegalite}), we get that
\begin{lemm}\label{minimal}
    For all diagram $D$, the matrix $D-\Gamma(D)$ has non negative integer entries.
\end{lemm}
For $A$ and $B$ two $k\times n$ matrix of non negative integers, let us say that $A\leq B$ if and only if $B-A$ has non negative entries. This is a partial order on $k\times n$ matrices of non negative integers. Putting together our observations of section \ref{compact} with (\ref{idemp_gamma}) and Lemma \ref{minimal}, we get
\begin{prop}
     For each permutation $\sigma$ in $S_n$, there is a unique compact diagram, $D_\sigma$, in the class of diagrams classified by $\sigma$. Moreover, $D_\sigma$ is minimal. This is to say that we have the inequality
\begin{equation}\label{minimalite}
    D_\sigma \leq D,
\end{equation}
as $2\times n$ non negative integer matrices, for all $D$ classified by $\sigma$.  
\end{prop}

\begin{proof}
The only part that remains to be checked is that $\Gamma(D)$ is indeed classified by $\sigma$, but this readily follows from the definition of $\Gamma(D)$.
\end{proof}

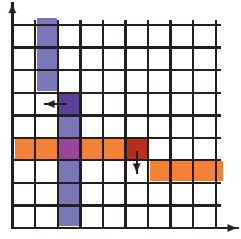
\begin{figure}[ht] 
\centering
\begin{picture}(11,11)(0,0)\thinlines
\put(0,0){\vector(1,0){10}}
\put(0,0){\vector(0,1){10}}
\put(2.07,5.5){\bleuf{\linethickness{2.5mm}\line(1,0){.85}}}
\put(2.51,0.1){\bleup{\linethickness{2.5mm}\line(0,1){4.8 }}}
\put(1.51,6.1){\bleup{\linethickness{2.5mm}\line(0,1){3.2 }}}
\put(2.4,5.5){\black{\vector(-1,0){1}}}
\put(5.07,3.5){\rouge{\linethickness{2.5mm}\line(1,0){.85}}}
\put(6.1,2.51){\rose{\linethickness{2.5mm}\line(1,0){3.2}}}
\put(0.13,3.51){\rose{\linethickness{2.5mm}\line(1,0){4.8}}}
\put(5.5,3.4){\black{\vector(0,-1){1}}}
\put(2.09,3.51){\violet{\linethickness{2.5mm}\line(1,0){.85}}}
\multiput(0,0)(0,1){10}{\line(1,0){9.2}}
\multiput(0,0)(1,0){10}{\line(0,1){9.2}}
\end{picture}
\caption{Constraint intervals for diagram compacting moves.}
\label{fig:intervalles vides pour mouv permis}
\end{figure}

Another approach to compactification of diagrams is through step by step transformations of diagrams that ultimately turns them into compact diagrams. Namely, for a cell $c=(a,b)$, we define the {\em left} translation of $c$
   $$c\ \leadsto\  \triangleleft (c):=(a-1,b),$$ 
and {\em down}  translation of $c$
 $$c\ \leadsto\  \triangledown (c):=(a,b-1).$$
A {\em compacting move} consists in replacing $c$, in a diagram $D$,  either by the cell $\triangleleft (c)$ or $\triangledown (c)$ when some conditions described below are fulfilled. We respectively denote  $\triangleleft_c(D)$ and $\triangledown_c(D)$ the resulting diagrams. To describe the constraints on compacting moves, we first associate to each cell for $c$ in $D$ some constraints intervals (see Figure \ref{fig:intervalles vides pour mouv permis}): 
    $$ {\rm Vert}(c,D):=\{ (a,b)\in D\ |\  \triangleleft(c)\prec  (a,b) 
                          \prec c\ \}$$           
and
     $${\rm Horiz}(c,D):= {\rm Vert}(c^{-1},D^{-1})^{-1},$$
where $(a,b)^{-1}=(b,a)$.

We can then set
\begin{equation}
\label{mouv_gauche}
   \triangleleft_c(D):= (D\setminus \{c\}) \cup \{\triangleleft (c)\}, 
\end{equation}
 if $\triangleleft (c)$ is in $\N\times \N$ and ${\rm Vert}(c,D)$ is empty. 
The reasoning behind this is that $\triangleleft_c$ moves the cell $c=(a,b)$ one unit to the left, if this change does not modify the underlying classifying permutation. In other words, the condition that ${\rm Vert}(c,D)$ be empty  insures that
\begin{equation}
    \sigma(\triangleleft_c(D))=\sigma(D).
 \end{equation}
  In a similar manner, we set
\begin{equation}
\label{mouv_enbas}
   \triangledown_c(D):=(D\setminus \{c\}) \cup \{\triangledown (c)\},
 \end{equation}
 if $\triangledown (c)$ is in $\N\times \N$ and  ${\rm Horiz}(c,D)$ is empty. To summarize, 
    $$ \triangleleft_c(D)\simeq D,\qquad {\rm and} \qquad \triangledown_c(D)\simeq D,$$
when the respective conditions are met. It is easy to check that a diagram $D$ is  compact if and only if there is no possible compacting move. 

\section{The fundamental bijection.}\label{bijection_fund}
The reason behind the introduction of compact diagrams is the following theorem\footnote{This theorem  is closely related, although in a more effective form, to a theorem of Garsia and Gessel (see \cite{garsia_gessel}) on bipartite partitions.} that will play a key role in showing that the set $\mathcal{M}_n$ is a basis of the trivial part of $R_{S_n\times S_n}$.
\begin{theo}\label{principal}
  There is a natural bijection between $n$-cell diagrams, and triplets  $(D_{\sigma},\lambda,\mu)$,
where $\sigma$ is the classifying permutation of $D$ and both $\lambda$ and $\mu$ are integer partitions having at most $n$ non zero parts. Moreover,
  $$\omega(D)=\omega(D_\sigma)+(|\lambda|,|\mu|),$$
 where $\omega(D):=\sum_{(a,b)\in D} (a,b)$.
 
 \end{theo}
\begin{proof} Let us denote $\varphi$ the bijection in question, and set 
\begin{equation}\label{triplet}
    \varphi(D):=(D_{\sigma},\lambda,\mu),
\end{equation}
with $\lambda$ and $\mu$ defined as follows. We simply set
  $$\lambda_{n+1-i} := a_i - d_i(\sigma),$$
and
  $$\mu_{n+1-\sigma(i)}:= b_i-d_{\sigma(i)}(\sigma^{-1}).$$
The particular indexing in these definitions insures that parts of $\lambda$ and $\mu$ are in decreasing order. 
Recall that
\begin{eqnarray*}
    D_\sigma&=&(\alpha(\sigma),\beta(\sigma))\\ 
                   &=& \{ (d_i(\sigma),d_{\sigma(i)}(\sigma^{-1})\ |\ 1\leq i\leq n\ \}.
\end{eqnarray*}
It follows that $\lambda$ and $\mu$ are respectively obtained by reordering in decreasing order the entries of the first and second line of the matrix $D-D_\sigma(i)$. This makes it evident that $\varphi$ is a bijection.
The compact diagram clearly corresponds to the case when $\lambda = \mu = 0$.
\end{proof}
\section{Basis of symmetric coinvariants.}\label{base coinvariants}
We are now in a position to prove that the set $\mathcal{M}_n$ (see  (\ref{la_base_triviale})) actually is a basis of the trivial isotypic component of
the space $R_{S_n\times S_n}$. This will simply follow from the following proposition.

\begin{prop}
  The set of diagonally symmetric polynomial 
  \begin{equation}\label{autre_base}
      \{\ m_\lambda(\mathbf{x}) m_\mu(\mathbf{y}) M_{(\alpha(\sigma),\beta(\sigma))}(\mathbf{x},\mathbf{y})\ |\ \ell(\lambda)\leq n,\ \ell(\mu)\leq n,\ {\rm and}\ \sigma\in S_n\ \}
  \end{equation}
\end{prop}
is linearly independent. 

\begin{proof}
 We need only observe that, for the lexicographic monomial order, the leading monomial of $m_\lambda(\mathbf{x}) m_\mu(\mathbf{y}) M_{(\alpha(\sigma),\beta(\sigma))}(\mathbf{x},\mathbf{y})$ is simply $X^\gamma Y^\delta$, where 
      $$\varphi(D)=(D_\sigma,\lambda,\mu),$$
 with $D=(\gamma,\delta)$.
\end{proof}

\section{Strict diagrams.}
\label{sect:compacts de diagramme a la gessel}
As we have seen in section \ref{SnxSn}, the function 
 \begin{equation}\label{flip_x}
    \downarrow_{\mathbf{x}} P(\mathbf{x},\mathbf{y})= P(\partial\mathbf{x}, \mathbf{y}) \Delta_n(\mathbf{x})
 \end{equation}
 establishes a sign twisting automorphism of the $S_n$-module of $S_n\times S_n$-harmonics. This is reflected in part by a natural bijection between the set of compact diagrams, and a natural indexing set for a basis of the alternating part of this same $S_n$-module.
This is also equivalent to a description of an explicit basis for the alternating $S_n$-isotypic component of the coinvariant module of $S_n\times S_n$.
 
Once again, we are looking for a family of special  (compact) ``strict''  diagrams indexed by permutations, together with a bijective encoding of general strict diagrams as triples $(D^{\rm s}_\sigma,\lambda,\mu)$, with $D^{\rm s}_\sigma$ one of these special strict diagrams, and $\lambda$ and $\mu$ partitions having at most $n$ parts.  These diagrams are simply obtained by translating, in term of diagrams, the effect of the sign twisting automorphism in (\ref{flip_x}). This takes the form of the following transform on compact diagrams:
 \begin{equation}\label{flip_diagram}
      D_\sigma\ \mapsto\ D_\sigma^{\rm s}:=(0,\sigma-1^n)+ {\rm J} D_\sigma,
         \qquad {\rm with}\quad  J:=\begin{pmatrix} 1&0\\0&-1
 \end{pmatrix}
 \end{equation}
 where $\sigma-1^n:=(\sigma(1)-1,\sigma(2)-1,\ldots,\sigma(n)-1)$. One easily verifies easily that the resulting diagrams have all their cells distinct, and it is clear that
     $$\sum_{(a,b)\in D^s_{\sigma}} (a,b)= \left({\rm maj}(\sigma),{\textstyle{n\choose 2}}-{\rm maj}(\sigma^{-1})\right).$$
The following table illustrates the result of this process.
{\footnotesize  \begin{eqnarray*}
     123 & \begin{pmatrix} 0&0&0 \\ 0&0&0\end{pmatrix} &\mapsto
               \begin{pmatrix} 0&0&0 \\ 0&1&2\end{pmatrix} \\
     132 & \begin{pmatrix} 0&0&1 \\ 0&1&0\end{pmatrix} &\mapsto
               \begin{pmatrix} 0&0&1 \\ 0&1&1\end{pmatrix} \\ 
     213 & \begin{pmatrix} 0&1&1 \\ 1&0&1\end{pmatrix} &\mapsto
               \begin{pmatrix} 0&1&1 \\ 0&0&1\end{pmatrix} \\ 
     231 & \begin{pmatrix} 0&0&1 \\ 1&1&0\end{pmatrix} &\mapsto
               \begin{pmatrix} 0&0&1 \\ 0&1&0\end{pmatrix} \\
     312 & \begin{pmatrix} 0&1&1 \\ 1&0&0\end{pmatrix} &\mapsto
               \begin{pmatrix} 0&1&1 \\ 1&0&1\end{pmatrix} \\ 
     321 & \begin{pmatrix} 0&1&2 \\ 2&1&0\end{pmatrix} &\mapsto
               \begin{pmatrix} 0&1&2 \\ 0&0&0\end{pmatrix} \\ 
 \end{eqnarray*}}
Strictly compact diagrams, for all permutations in  $S_4$, are depicted in Figure \ref{fig:diagrams pour 4}.

{\tiny \begin{figure}[ht] 
$$\begin{matrix}\qquad
\begin{picture}(1,4)(0,0)
\put(0.05,3.5){\carre}
\put(0.05,2.5){\carre}
\put(0.05,1.5){\carre}
\put(0.05,0.5){\carre}
\end{picture} \qquad & \qquad
\begin{picture}(2,3)(0,0)
\put(0.05,0.5){\carre}
\put(0.05,2.5){\carre}
\put(0.05,1.5){\carre}
\put(1.05,2.5){\carre}
\end{picture}  \qquad & \qquad
\begin{picture}(2,3)(0,0)
\put(0.05,0.5){\carre}
\put(0.05,1.5){\carre}
\put(1.05,1.5){\carre}
\put(1.05,2.5){\carre}
\end{picture}  \qquad & \qquad
\begin{picture}(2,3)(0,0)
\put(0.05,0.5){\carre}
\put(0.05,2.5){\carre}
\put(0.05,1.5){\carre}
\put(1.05,1.5){\carre}
\end{picture} \qquad \\
  1234 & 1243 & 1324 & 1342 \\ 
\begin{picture}(2,3.5)(0,0)
\put(0.05,2.5){\carre}
\put(1.05,1.5){\carre}
\put(1.05,2.5){\carre}
\put(0.05,0.5){\carre}
\end{picture} &
\begin{picture}(3,2)(0,0)
\put(0.05,0.5){\carre}
\put(0.05,1.5){\carre}
\put(1.05,1.5){\carre}
\put(2.05,1.5){\carre}
\end{picture} &
\begin{picture}(2,3)(0,0)
\put(0.05,0.5){\carre}
\put(1.05,1.5){\carre}
\put(1.05,0.5){\carre}
\put(1.05,2.5){\carre}
\end{picture} &
\begin{picture}(3,2)(0,0)
\put(0.05,1.5){\carre}
\put(1.05,1.5){\carre}
\put(1.05,0.5){\carre}
\put(2.05,0.5){\carre}
\end{picture} \\
1423 & 1432 & 2134 & 2143 \\ 
\begin{picture}(2,3.5)(0,0)
\put(1.05,0.5){\carre}
\put(0.05,0.5){\carre}
\put(0.05,1.5){\carre}
\put(1.05,2.5){\carre}
\end{picture} &
\begin{picture}(2,3)(0,0)
\put(0.05,0.5){\carre}
\put(0.05,2.5){\carre}
\put(0.05,1.5){\carre}
\put(1.05,0.5){\carre}
\end{picture} &
\begin{picture}(2,2)(0,0)
\put(0.05,1.5){\carre}
\put(0.05,0.5){\carre}
\put(1.05,1.5){\carre}
\put(1.05,0.5){\carre}
\end{picture} &
\begin{picture}(3,2)(0,0)
\put(0.05,0.5){\carre}
\put(0.05,1.5){\carre}
\put(1.05,1.5){\carre}
\put(2.05,0.5){\carre}
\end{picture} \\
2314 & 2341 & 2413 & 2431 \\ 
\begin{picture}(2,3.5)(0,0)
\put(1.05,2.5){\carre}
\put(1.05,1.5){\carre}
\put(1.05,0.5){\carre}
\put(0.05,1.5){\carre}
\end{picture} &
\begin{picture}(3,3)(0,0)
\put(0.05,1.5){\carre}
\put(1.05,2.5){\carre}
\put(1.05,0.5){\carre}
\put(2.05,1.5){\carre}
\end{picture} &
\begin{picture}(3,2)(0,0)
\put(0.05,0.5){\carre}
\put(1.05,0.5){\carre}
\put(2.05,0.5){\carre}
\put(2.05,1.5){\carre}
\end{picture} &
\begin{picture}(3,2)(0,0)
\put(0.05,0.5){\carre}
\put(1.05,1.5){\carre}
\put(1.05,0.5){\carre}
\put(2.05,0.5){\carre}
\end{picture} \\
3124 & 3142 & 3214 & 3241 \\ 
\begin{picture}(2,3.5)(0,0)
\put(0.05,2.5){\carre}
\put(0.05,1.5){\carre}
\put(1.05,0.5){\carre}
\put(1.05,1.5){\carre}
\end{picture} &
\begin{picture}(3,2)(0,0)
\put(0.05,0.5){\carre}
\put(1.05,0.5){\carre}
\put(2.05,0.5){\carre}
\put(0.05,1.5){\carre}
\end{picture} &
\begin{picture}(2,3)(0,0)
\put(0.05,2.5){\carre}
\put(1.05,1.5){\carre}
\put(1.05,0.5){\carre}
\put(1.05,2.5){\carre}
\end{picture} &
\begin{picture}(3,2)(0,0)
\put(0.05,1.5){\carre}
\put(1.05,1.5){\carre}
\put(1.05,0.5){\carre}
\put(2.05,1.5){\carre}
\end{picture} \\
3412 & 3421 & 4123 & 4132 \\ 
\begin{picture}(3,2.5)(0,0)
\put(0.05,1.5){\carre}
\put(1.05,0.5){\carre}
\put(2.05,0.5){\carre}
\put(2.05,1.5){\carre}
\end{picture} &
\begin{picture}(3,2)(0,0)
\put(0.05,1.5){\carre}
\put(1.05,1.5){\carre}
\put(1.05,0.5){\carre}
\put(2.05,0.5){\carre}
\end{picture} &
\begin{picture}(3,2)(0,0)
\put(0.05,1.5){\carre}
\put(2.05,0.5){\carre}
\put(1.05,1.5){\carre}
\put(2.05,1.5){\carre}
\end{picture} &
\begin{picture}(4,1)(0,0)
\put(0.05,0.5){\carre}
\put(1.05,0.5){\carre}
\put(2.05,0.5){\carre}
\put(3.05,0.5){\carre}
\end{picture}\\
4213 & 4231 & 4312 & 4321
\end{matrix}$$
\caption{Strictly compact diagrams, $n=4$.}
\label{fig:diagrams pour 4}
\end{figure}
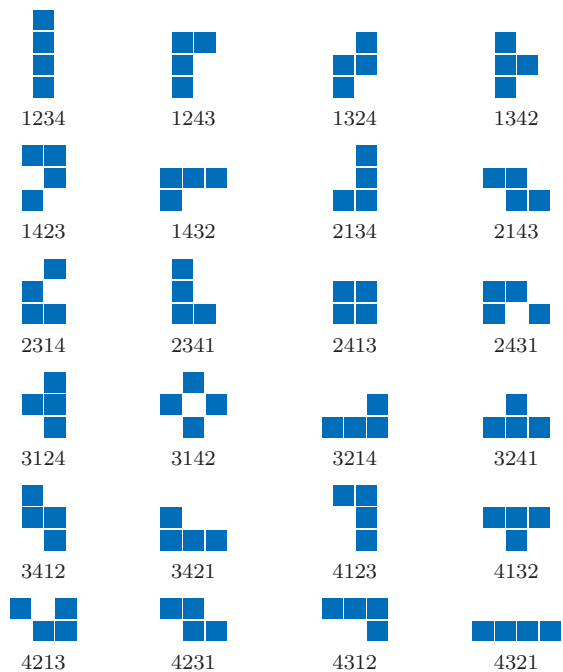}

\begin{figure}[ht] 
\centering
\begin{picture}(11,11)(0,0)\thinlines
\put(0,0){\vector(1,0){10}}
\put(0,0){\vector(0,1){10}}
\put(2.07,5.5){\bleuf{\linethickness{2.5mm}\line(1,0){.85}}}
\put(2.51,0.1){\bleup{\linethickness{2.5mm}\line(0,1){4.8 }}}
\put(1.51,6.1){\bleup{\linethickness{2.5mm}\line(0,1){3.2 }}}
\put(2.4,5.5){\black{\vector(-1,0){1}}}
\put(5.07,3.5){\rouge{\linethickness{2.5mm}\line(1,0){.85}}}
\put(6.1,3.51){\rose{\linethickness{2.5mm}\line(1,0){3.2}}}
\put(0.13,2.51){\rose{\linethickness{2.5mm}\line(1,0){5.8}}}
\put(5.5,3.4){\black{\vector(0,-1){1}}}
\put(2.09,2.51){\violet{\linethickness{2.5mm}\line(1,0){.85}}}
\multiput(0,0)(0,1){10}{\line(1,0){9.2}}
\multiput(0,0)(1,0){10}{\line(0,1){9.2}}
\end{picture}
\caption{Constraint intervals for strict compacting moves.}
\label{strict_constraints}
\end{figure}
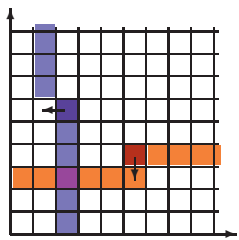

Just as in our previous case, we can describe strictly compact diagrams directly in term of their indexing permutation, or in term of {\em strict compacting} moves. These are the same moves as before, but the constraint intervals are modified as follows. We now set
    $$ {\rm Vert}_s(c,D):=\{ (a,b)\in D\ |\  \triangleleft(c)\preceq (a,b)
                          \prec c\ \}$$     
and
    $$ {\rm Horiz}_s(c,D):=\{ (a,b)\in D\ |\  (b,a)\preceq \triangleleft(c^{-1}),\ {\rm or}\ 
                          c^{-1} \prec (b,a)\ \}$$           
These two intervals considered are illustrated in
Figure~\ref{strict_constraints}.

Once again, there is a natural bijection between the set of strict $n$-cells diagrams, and the
set of triples $(D^{\rm s}_\sigma,\lambda,\mu)$, with $\lambda$ and $\mu$ partitions having at most $n$ parts.  Many nice properties of strict diagrams that have been explored in  \cite{lamontagne}.

\section{Final remarks.}\label{final}
The constructions and results for diagrams and strict diagrams afford a common generalization that allows a combinatorial description of each isotypic component of the $S_n$-module of $S_n\times S_n$-coinvariants.  More precisely, there is a notion of compact diagrams indexed by pairs $(\sigma,\tau)$, with $\sigma$ a permutation and $\tau$ a standard tableau, where the charge statistic plays a natural role. There are also similar results for the decomposition of the diagonal action of $S_n$ on  $S_n^k$-coinvariants, as well as for other Coxeter groups such as $B_n$ and $D_n$. All these generalization will be the subject of an upcoming paper.


\end{document}